\documentclass{amsart}
\usepackage{amsfonts}
\usepackage{amssymb}
\usepackage{amsmath}

\newtheorem{theorem}{Theorem}[section]
\theoremstyle{plain}

\newtheorem{corollary}[theorem]{Corollary}

\numberwithin{equation}{section}
\input{tcilatex}

\begin{document}

\begin{center}
{\Large Univalence criterion for meromorphic functions and Loewner chains }

\textbf{Erhan DEN\.{I}Z and Halit ORHAN}

\textit{Dedicated to Professor H. M. Srivastava on the Occasion of his
Seventieth Birth Anniversary}
\end{center}

\textbf{Abstract.} The object of the present paper is to obtain a more
general condition for univalence of meromorphic functions in the $\mathbb{U}%
^{\ast }$. The significant relationships and relevance with other results
are also given. A number of known univalent conditions would follow upon
specializing the parameters involved in our main results.

\begin{flushleft}

\end{flushleft}

\section{Introduction}

We denote by $\mathbb{U}_{r}$ the disk $\left\{ z\in {\mathbb{C}}%
:\;\left\vert z\right\vert <r\right\} ,$ where $0<r\leq 1$, by $\mathbb{U}=%
\mathbb{U}_{1}$ the open unit disk of the complex plane and $\mathbb{U}%
^{\ast }=%
\mathbb{C}
\diagdown \overline{\mathbb{U}},$ where $\overline{\mathbb{U}}$ is closure
of $\mathbb{U}.$

Let $\mathcal{A}$ denote the class of all analytic functions in the open
unit disk $\mathbb{U}$ normalized by%
\begin{equation*}
f(z)=z+a_{2}z^{2}+...\text{ \ \ }\left( z\in \mathbb{U}\right) 
\end{equation*}%
and we denote by $\mathcal{S}$ the subclass of $\mathcal{A}$ consisting of
functions which are also univalent in $\mathbb{U}.$ Closely related to $%
\mathcal{S}$ is the class $\sum $ of all meromorphic functions in $\mathbb{U}%
^{\ast }$ by%
\begin{equation*}
f(\zeta )=b\zeta +b_{0}+\frac{b_{1}}{\zeta }+...\text{ \ \ }\left( \zeta \in 
\mathbb{U}^{\ast }\right) 
\end{equation*}%
and $\sum_{0}$ stands for all functions from $\sum $ with normalization $b=1$
and $b_{0}=0.$ These classes have been one of the important subjects of
research in complex analysis especially, Geometric Function Theory for a
long time (see, for details, \cite{Sr}).

Two of the most important and known univalence criteria for analytic
functions defined in $\mathbb{U}^{\ast }$ were obtained by Becker \cite{Be}
and Nehari \cite{Ne}. Some extensions of these two criteria were given by
Lewandowski \cite{Le1}, \cite{Le2} and Ruscheweyh \cite{Ru}. During the
time, unlike there were obtained a lot of univalence criteria by Miazga and
Wesolowski \cite{MiWe}, Wesolowski \cite{We}, Kanas and Srivastava \cite%
{KaSr} and Deniz and Orhan \cite{DeOr}.

In the present paper we consider a general univalence criterion for
functions $f$ belonging to the class $\sum $ in terms of the Schwarz
derivative defined by%
\begin{equation*}
{S_{f}(z)=}\left( \frac{f^{\prime \prime }(z)}{f^{\prime }(z)}\right)
^{\prime }-\frac{1}{2}\left( \frac{f^{\prime \prime }(z)}{f^{\prime }(z)}%
\right) ^{2}.
\end{equation*}

\section{Loewner chains and related theorem}

Before proving our main theorem we need a brief summary of the method of
Loewner chains.

Let $\mathcal{L}(z,t)=a_{1}(t)z+a_{2}(t)z^{2}+...,$ $a_{1}(t)\neq 0,$ be a
function defined on $\mathbb{U}\times \lbrack 0,\infty )$, where $a_{1}(t)$
is a complex-valued, locally absolutely continuous function on $[0,\infty ).$
$\mathcal{L}(z,t)$ is called a Loewner chain if $\mathcal{L}(z,t)$ satisfies
the following conditions;

\begin{enumerate}
\item[(i)] $\mathcal{L}(z,t)$ is analytic and univalent in $\mathbb{U}$ for
all $t\in \lbrack 0,\infty )$

\item[(ii)] $\mathcal{L}(z,t)\prec \mathcal{L}(z,s)$ for all $0\leq t\leq
s<\infty $,
\end{enumerate}

where the symbol \textquotedblright $\prec $ \textquotedblright\ stands for
subordination. If $a_{1}(t)=e^{t}$ then we say that $\mathcal{L}(z,t)$ is a 
\textit{standard Loewner chain}.

In order to prove our main results we need the following theorem due to
Pommerenke \cite{Po1} (also see \cite{Po2}). This theorem is often used to
find out univalency for an analytic function, apart from the theory of
Loewner chains;

\begin{theorem}
\label{t0}\textbf{\ }\textit{Let }$\mathcal{L}%
(z,t)=a_{1}(t)z+a_{2}(t)z^{2}+...$\textit{\ be analytic in }$\mathbb{U}_{r}$%
\textit{\ for all }$t\in \lbrack 0,\infty )$\textit{$.$ S}uppose that;

\begin{enumerate}
\item[(i)] \textit{\ }$\mathcal{L}(z,t)$ is a \textit{locally absolutely
continuous function in the interval }$[0,\infty ),$\textit{\ and locally
uniformly with respect to }$\mathbb{U}_{r}.$\textit{\ }

\item[(ii)] $a_{1}(t)$ is a complex valued continuous function on $[0,\infty
)$ such that $a_{1}(t)\neq 0,$ $\left\vert {a_{1}(t)}\right\vert \rightarrow
\infty $\textit{\ for }$t\rightarrow \infty $ and%
\begin{equation*}
\left\{ \frac{{\mathcal{L}(z,t)}}{a_{1}(t)}\right\} _{t\in \lbrack 0,\infty
)}
\end{equation*}%
\textit{\ forms a normal family of functions in }$\mathbb{U}_{r}.$

\item[(iii)] There exists an analytic function $p:\mathbb{U}\times \lbrack
0,\infty )\rightarrow 
\mathbb{C}
$ satisfying $\func{Re}{p(z,t)}>0$\textit{\ for all }$z\in \mathbb{U},\;t\in
\lbrack 0,\infty )$ and 
\begin{equation}
z\frac{\partial \mathcal{L}(z,t)}{\partial z}=p(z,t)\frac{\partial \mathcal{L%
}(z,t)}{\partial t},\quad z\in \mathbb{U}_{r},\text{ }t\in \lbrack 0,\infty
).  \label{1}
\end{equation}%
\textit{Then, for each }$t\in \lbrack 0,\infty ),$\textit{\ the function }$%
\mathcal{L}(z,t)$\textit{\ has an analytic and univalent extension to the
whole disk }$\mathbb{U}$ or the function $\mathcal{L}(z,t)$ is a Loewner
chain.
\end{enumerate}
\end{theorem}

The equation (\ref{1}) is called the generalized Loewner differential
equation.

\section{Univalence criterion for the functions belonging to the class $\sum 
$}

In this section, making use of the Theorem \ref{t0}, we obtain an univalence
criterion for meromorphic functions. The method of prove is based on Theorem %
\ref{t0} and on construction of a suitable Loewner chain.

\begin{theorem}
\label{t1}\textbf{\ }Let \textbf{\ }$f,g\in \sum $ be locally univalent
functions in $\mathbb{U}^{\ast }.$ If there exists an analytic function $h$
such that $\func{Re}h(\zeta )\geq \frac{1}{2}$ and $h(\zeta )=1+\frac{h_{2}}{%
\zeta ^{2}}+...$ for$\ \ \zeta \in \mathbb{U}^{\ast },$ and for arbitrary $%
\alpha \in 
\mathbb{C}
$ we have%
\begin{eqnarray}
&&\left\vert {\frac{1-h(\zeta )}{h(\zeta )}\left\vert \zeta \right\vert
^{2}-(\left\vert \zeta \right\vert ^{2}-1)\left[ {\frac{\zeta {h}^{\prime
}(\zeta )}{h(\zeta )}+(1-2\alpha )\frac{\zeta {f}^{\prime \prime }(\zeta )}{{%
f}^{\prime }(\zeta )}+2\alpha \frac{\zeta {g}^{\prime \prime }(\zeta )}{{g}%
^{\prime }(\zeta )}}\right] }\right.  \label{3.1} \\
&&\left. {\ +\alpha (\left\vert \zeta \right\vert ^{2}-1)^{2}\frac{\zeta }{%
\overline{\zeta }}h(\zeta )\left[ {\left( {\alpha -\frac{1}{2}}\right)
\left( {\frac{{f}^{\prime \prime }(\zeta )}{{f}^{\prime }(\zeta )}-\frac{{g}%
^{\prime \prime }(\zeta )}{{g}^{\prime }(\zeta )}}\right) ^{2}+{S_{f}(\zeta
)-S_{g}(\zeta )}}\right] }\right\vert \leqslant 1  \notag
\end{eqnarray}%
for all$\ \zeta \in \mathbb{U}^{\ast },$ then $f$ is univalent in $\mathbb{U}%
^{\ast }.$
\end{theorem}

\begin{proof}
Without loss of generality we can consider the functions of the form%
\begin{equation*}
f(\zeta )=\zeta +\frac{{a_{1}}}{\zeta }+...\text{ and \ }g(\zeta )=\zeta +%
\frac{{b_{1}}}{\zeta }+...
\end{equation*}%
\noindent since the Schwarzian derivative is invariant under M\"{o}bius
transformations. Consider the functions defined by%
\begin{equation}
v(\zeta )=\left[ {\frac{{g}^{\prime }(\zeta )}{{f}^{\prime }(\zeta )}}\right]
^{\alpha }=1+\frac{v_{2}}{\zeta ^{2}}+...,\quad \alpha \in \mathbb{C}
\label{3.2}
\end{equation}%
\noindent where we choose this branch of the power $(\cdot )^{\alpha },$
which for $\zeta \rightarrow \infty $ has value $1,$ and%
\begin{equation}
u(\zeta )=f(\zeta )v(\zeta )=\zeta +\frac{u_{2}}{\zeta }+....  \label{3.3}
\end{equation}%
The functions $u$ and $v$ are meromorphic in $U^{\ast }$ since $f$ and $g$
do not have multiple poles and ${f}^{\prime }$ and ${g}^{\prime }$ are
different from zero.

For all $t\in \left[ {0,\infty }\right) $ and $\frac{1}{\zeta }=z\in \mathbb{%
U}$ the function $f:\mathbb{U}_{r}\times \left[ {0,\infty }\right)
\rightarrow \mathbb{C}$ defined formally by%
\begin{eqnarray}
f(z,t) &=&\left[ {\frac{u\left( {\frac{e^{t}}{z}}\right) +(e^{-t}-e^{t})%
\frac{1}{z}h\left( {\frac{e^{t}}{z}}\right) {u}^{\prime }\left( {\frac{e^{t}%
}{z}}\right) }{v\left( {\frac{e^{t}}{z}}\right) +(e^{-t}-e^{t})\frac{1}{z}%
h\left( {\frac{e^{t}}{z}}\right) {v}^{\prime }\left( {\frac{e^{t}}{z}}%
\right) }}\right] ^{-1}  \label{3.4} \\
&=&e^{t}z+\Psi (e^{-pt},z^{2}),\quad p=1,2,...  \notag
\end{eqnarray}%
\noindent is analytic in $\mathbb{U}$ since $\Psi (e^{-pt},z^{2})$ is
analytic function in $\mathbb{U}$ for each fixed $t\in \left[ {0,\infty }%
\right) $ and $p=1,2,....$ From (\ref{3.4}) we have $a_{1}(t)=e^{t}$ and%
\begin{equation*}
\underset{t\rightarrow \infty }{\lim }\left\vert {a_{1}(t)}\right\vert =%
\underset{t\rightarrow \infty }{\lim }e^{t}=\infty .
\end{equation*}%
After simple calculation we obtain, for each $z\in \mathbb{U},$%
\begin{equation*}
\underset{t\rightarrow \infty }{\lim }\frac{f(z,t)}{e^{t}}=\underset{%
t\rightarrow \infty }{\lim }\left\{ {z+\Psi (e^{-(p+1)t},z^{2})}\right\} =z.
\end{equation*}%
The limit function $k(z)=z$ belongs to the family $\left\{ f(z,t)\diagup
e^{t}:\;t\in \left[ {0,\infty }\right) \right\} ;$ then, there exists a
number $r_{0}$ $(0<r_{0}<1)$ that in every closed disk $\mathbb{U}_{r_{0}},$
there exists a constant $K_{0}>0,$ such that%
\begin{equation*}
\left\vert {\frac{f(z,t)}{e^{t}}}\right\vert <K_{0},\quad z\in \mathbb{U}%
_{r_{0}},\;t\in \left[ {0,\infty }\right)
\end{equation*}%
uniformly in this disk, provided that $t$ is sufficiently large. Thus, by
Montel's Theorem, $\left\{ f(z,t)\diagup e^{t}\right\} $ forms a normal
family\ in each disk $\mathbb{U}_{r_{0}}.$

Since the function $\Psi (e^{-pt},z^{2})$ is analytic in $\mathbb{U},$ $\Psi
^{(k)}(e^{-pt},z^{2})$ $k\in \mathbb{N}_{0}=\{0,1,2...\}$ is \ continuous on
the compact set, so $\Psi ^{(k)}(e^{-pt},z^{2}),$ $k\in \mathbb{N}_{0}$ is
bounded function. Thus for all fixed $T>0$, we can write ${e^{t}}<e^{T}$ and
we obtain that for all fixed numbers $t\in \left[ {0,T}\right] \subset \left[
{0,\infty }\right) ,$ there exists a constant $K_{1}>0$ such that 
\begin{equation*}
\left\vert {\frac{\partial f(z,t)}{\partial t}}\right\vert <K_{1},\text{ \ \ 
}\forall z\in \mathbb{U}_{r_{0}},\;t\in \left[ {0,T}\right] .
\end{equation*}%
\noindent Therefore, the function $f(z,t)$ is locally absolutely continuous
in $\left[ {0,\infty }\right) $; locally uniformly with respect to $\mathbb{U%
}_{r_{0}}.$

After simple calculations from (\ref{3.4}) we obtain%
\begin{eqnarray}
&&\frac{\partial f(z,t)}{\partial z}  \label{3.41} \\
&=&\frac{1}{z}\frac{e^{t}}{z}\left\{ {\left( {1+(e^{-2t}-1)\left[ {h\left( 
\frac{{{e^{t}}}}{{z}}\right) +\frac{e^{t}}{z}{h}^{\prime }\left( \frac{{{%
e^{t}}}}{{z}}\right) }\right] }\right) }\left( {{u}^{\prime }v-{v}^{\prime }u%
}\right) \right.  \notag \\
&&\left. {\ +(e^{-2t}-1)\frac{e^{t}}{z}h\left( \frac{{{e^{t}}}}{{z}}\right)
\left( {{u}^{\prime \prime }v-{v}^{\prime \prime }u}\right) +(e^{-2t}-1)^{2}%
\frac{e^{2t}}{z^{2}}h^{2}\left( \frac{{{e^{t}}}}{{z}}\right) \left( {{u}%
^{\prime \prime }{v}^{\prime }-{v}^{\prime \prime }{u}^{\prime }}\right) }%
\right\}  \notag \\
&&\times f^{2}(z,t)/\left[ {v\left( {\frac{e^{t}}{z}}\right) +(e^{-t}-e^{t})%
\frac{1}{z}h\left( {\frac{e^{t}}{z}}\right) {v}^{\prime }\left( {\frac{e^{t}%
}{z}}\right) }\right] ^{2}  \notag
\end{eqnarray}%
and%
\begin{eqnarray}
&&\frac{\partial f(z,t)}{\partial t}  \label{3.42} \\
&=&-\frac{e^{t}}{z}\left\{ {\left( {1-(e^{-2t}+1)h\left( \frac{{{e^{t}}}}{{z}%
}\right) +(e^{-2t}+1)\frac{e^{t}}{z}{h}^{\prime }\left( \frac{{{e^{t}}}}{{z}}%
\right) }\right) }\left( {{u}^{\prime }v-{v}^{\prime }u}\right) \right. 
\notag \\
&&\left. {\ +(e^{-2t}-1)\frac{e^{t}}{z}h\left( \frac{{{e^{t}}}}{{z}}\right)
\left( {{u}^{\prime \prime }v-{v}^{\prime \prime }u}\right) +(e^{-2t}-1)^{2}%
\frac{e^{2t}}{z^{2}}h^{2}\left( \frac{{{e^{t}}}}{{z}}\right) \left( {{u}%
^{\prime \prime }{v}^{\prime }-{v}^{\prime \prime }{u}^{\prime }}\right) }%
\right\}  \notag \\
&&\times f^{2}(z,t)/\left[ {v\left( {\frac{e^{t}}{z}}\right) +(e^{-t}-e^{t})%
\frac{1}{z}h\left( {\frac{e^{t}}{z}}\right) {v}^{\prime }\left( {\frac{e^{t}%
}{z}}\right) }\right] ^{2}  \notag
\end{eqnarray}%
where%
\begin{equation}
{u}^{\prime }v-{v}^{\prime }u={f}^{\prime }\left( {\frac{{g}^{\prime }}{{f}%
^{\prime }}}\right) ^{2\alpha },\text{ \ \ }\alpha \in 
\mathbb{C}
\label{3.5}
\end{equation}%
\begin{equation}
{u}^{\prime \prime }v-{v}^{\prime \prime }u=(1-2\alpha ){f}^{\prime \prime
}\left( {\frac{{g}^{\prime }}{{f}^{\prime }}}\right) ^{2\alpha }+2\alpha {g}%
^{\prime \prime }\left( {\frac{{g}^{\prime }}{{f}^{\prime }}}\right)
^{2\alpha -1},\text{ \ \ \ }\alpha \in 
\mathbb{C}
\label{3.6}
\end{equation}%
\begin{equation}
{u}^{\prime \prime }{v}^{\prime }-{v}^{\prime \prime }{u}^{\prime }=\alpha {f%
}^{\prime }\left( {\frac{{g}^{\prime }}{{f}^{\prime }}}\right) ^{2\alpha
}\left\{ {\left( {S_{f}-S_{g}}\right) +\left( {\alpha -\frac{1}{2}}\right)
\left( {\frac{{f}^{\prime \prime }}{{f}^{\prime }}-\frac{{g}^{\prime \prime }%
}{{g}^{\prime }}}\right) }\right\} ,\text{ \ \ }\alpha \in 
\mathbb{C}
\label{3.7}
\end{equation}%
and $u,v,u^{\prime },v^{\prime },u^{\prime \prime },v^{\prime \prime }$ are
calculated at $\frac{{{e^{t}}}}{{z}}.$

Consider the function $p:\mathbb{U}_{r}\times \left[ {0,\infty }\right)
\rightarrow \mathbb{C}$ for$\;0<r<r_{0}$ and $t\in \left[ {0,\infty }\right)
,$ defined by%
\begin{equation*}
p(z,t)={z\frac{\partial f(z,t)}{\partial z}}\diagup \frac{\partial f(z,t)}{%
\partial t}.
\end{equation*}%
From (\ref{3.41}) to (\ref{3.7}), we can easily see that the function $%
p(z,t) $ is analytic in $\mathbb{U}_{r},$ $0<r<r_{0}.$ If the function%
\begin{equation}
w(z,t)=\frac{p(z,t)-1}{p(z,t)+1}=\frac{\frac{z\partial f(z,t)}{\partial z}-%
\frac{\partial f(z,t)}{\partial t}}{\frac{z\partial f(z,t)}{\partial z}+%
\frac{\partial f(z,t)}{\partial t}}  \label{3.8}
\end{equation}%
\noindent is analytic in $\mathbb{U}\times \left[ {0,\infty }\right) $ and $%
\left\vert {w(z,t)}\right\vert <1,$ for all $z\in \mathbb{U}\;$and $t\in %
\left[ {0,\infty }\right) ,$ then $p(z,t)$ has an analytic extension with
positive real part $(\func{Re}p(z,t)>0)$ in $\mathbb{U},$ for all $t\in %
\left[ {0,\infty }\right) .$

To show this we write (\ref{3.41}) and (\ref{3.42}) in the equation (\ref%
{3.8}), then we obtain 
\begin{eqnarray}
&&w(z,t)  \label{3.9} \\
&=&\frac{2\frac{e^{t}}{z}\left\{ {\left( {1-h\left( \frac{{{e^{t}}}}{{z}}%
\right) +(e^{-2t}-1)\frac{e^{t}}{z}{h}^{\prime }\left( \frac{{{e^{t}}}}{{z}}%
\right) }\right) }\left( {{u}^{\prime }v-{v}^{\prime }u}\right) \right. }{%
2e^{-2t}\frac{e^{t}}{z}h\left( \frac{{{e^{t}}}}{{z}}\right) \left( {{u}%
^{\prime }v-{v}^{\prime }u}\right) }  \notag \\
&&+\frac{\left. {(e^{-2t}-1)\frac{e^{t}}{z}h\left( \frac{{{e^{t}}}}{{z}}%
\right) \left( {{u}^{\prime \prime }v-{v}^{\prime \prime }u}\right)
+(e^{-2t}-1)^{2}\frac{e^{2t}}{z^{2}}h^{2}\left( \frac{{{e^{t}}}}{{z}}\right)
\left( {{u}^{\prime \prime }{v}^{\prime }-{v}^{\prime \prime }{u}^{\prime }}%
\right) }\right\} }{2e^{-2t}\frac{e^{t}}{z}h\left( \frac{{{e^{t}}}}{{z}}%
\right) \left( {{u}^{\prime }v-{v}^{\prime }u}\right) }  \notag \\
&=&e^{2t}\left( {\frac{1-h\left( \frac{{{e^{t}}}}{{z}}\right) }{h\left( 
\frac{{{e^{t}}}}{{z}}\right) }}\right) +(1-e^{2t})\frac{e^{t}}{z}\left( {%
\frac{{h}^{\prime }\left( \frac{{{e^{t}}}}{{z}}\right) }{h\left( \frac{{{%
e^{t}}}}{{z}}\right) }+\frac{{u}^{\prime \prime }v-{v}^{\prime \prime }u}{{u}%
^{\prime }v-{v}^{\prime }u}}\right)  \notag \\
&&+e^{2t}(e^{-2t}-1)^{2}\frac{e^{2t}}{z^{2}}h\left( \frac{{{e^{t}}}}{{z}}%
\right) \frac{{u}^{\prime \prime }{v}^{\prime }-{v}^{\prime \prime }{u}%
^{\prime }}{{u}^{\prime }v-{v}^{\prime }u}  \notag
\end{eqnarray}%
and from (\ref{3.5})-(\ref{3.7}) for all $z\in \mathbb{U}$ and $t\in \left[ {%
0,\infty }\right) $%
\begin{eqnarray}
&&w(z,t)  \label{3.10} \\
&=&e^{2t}\left( {\frac{1-h\left( \frac{e^{t}}{z}\right) }{h\left( \frac{e^{t}%
}{z}\right) }}\right) +(1-e^{2t})\frac{e^{t}}{z}\left( {\frac{{h}^{\prime
}\left( \frac{e^{t}}{z}\right) }{h\left( \frac{e^{t}}{z}\right) }+(1-2\alpha
)\frac{{f}^{\prime \prime }\left( \frac{e^{t}}{z}\right) }{{f}^{\prime
}\left( \frac{e^{t}}{z}\right) }+2\alpha \frac{{g}^{\prime \prime }\left( 
\frac{e^{t}}{z}\right) }{{g}^{\prime }\left( \frac{e^{t}}{z}\right) }}\right)
\notag \\
&&+\alpha e^{2t}(e^{-2t}-1)^{2}\frac{e^{2t}}{z^{2}}h\left( \frac{e^{t}}{z}%
\right) \left( {\left( {S_{f}(\frac{e^{t}}{z})-S_{g}(\frac{e^{t}}{z})}%
\right) +\left( {\alpha -\frac{1}{2}}\right) \left( {\frac{{f}^{\prime
\prime }(\frac{e^{t}}{z})}{{f}^{\prime }(\frac{e^{t}}{z})}-\frac{{g}^{\prime
\prime }(\frac{e^{t}}{z})}{{g}^{\prime }(\frac{e^{t}}{z})}}\right) }\right) .
\notag
\end{eqnarray}%
The right hand side of the equation (\ref{3.10}) is equal to%
\begin{equation*}
w(z,0)=\frac{1-h\left( {\frac{1}{z}}\right) }{h\left( {\frac{1}{z}}\right) }
\end{equation*}%
for $t=0.$ Thus, from hypothesis of theorem for $\frac{1}{z}=\zeta \in 
\mathbb{U}^{\ast }$ we have%
\begin{equation*}
\left\vert {\frac{1-h\left( \zeta \right) }{h\left( \zeta \right) }}%
\right\vert \leqslant 1.
\end{equation*}%
Since $\left\vert {\frac{e^{t}}{z}}\right\vert \geqslant \left\vert {e^{t}}%
\right\vert >1$ for all $z\in \overline{\mathbb{U}}$ and $t>0,$ we find that 
$w(z,t)$ is an analytic function in $\overline{\mathbb{U}}.$ Then putting $%
\frac{e^{t}}{z}=\widetilde{\zeta }\in \mathbb{U}^{\ast }$, $\widetilde{\zeta 
}=\zeta e^{t}$, $\left\vert \widetilde{\zeta }\right\vert =e^{t}$ for $%
\left\vert z\right\vert =1,$ from (\ref{3.10}) by assumption (\ref{3.1})
replacing $\widetilde{\zeta }$ by $\zeta $ we have%
\begin{eqnarray*}
\left\vert {w(z,t)}\right\vert &=&\left\vert {\left\vert \zeta \right\vert
^{2}\left( {\frac{1-h\left( \zeta \right) }{h\left( \zeta \right) }}\right)
-(\left\vert \zeta \right\vert ^{2}-1)\left( {\frac{\zeta {h}^{\prime
}\left( \zeta \right) }{h\left( \zeta \right) }+(1-2\alpha )\frac{\zeta {f}%
^{\prime \prime }\left( \zeta \right) }{{f}^{\prime }\left( \zeta \right) }%
+2\alpha \frac{\zeta {g}^{\prime \prime }\left( \zeta \right) }{{g}^{\prime
}\left( \zeta \right) }}\right) }\right. \\
&&\left. {\;+\alpha (\left\vert \zeta \right\vert ^{2}-1)^{2}\frac{e^{2t}}{%
z^{2}}h\left( \zeta \right) \left( {\left( {S_{f}(\zeta )-S_{g}(\zeta )}%
\right) +\left( {\alpha -\frac{1}{2}}\right) \left( {\frac{{f}^{\prime
\prime }(\zeta )}{{f}^{\prime }(\zeta )}-\frac{{g}^{\prime \prime }(\zeta )}{%
{g}^{\prime }(\zeta )}}\right) }\right) }\right\vert \\
&\leqslant &1.
\end{eqnarray*}%
Therefore $\left\vert {w(z,t)}\right\vert <1$ for all $z\in \mathbb{U}\;$and 
$t\in \left[ {0,\infty }\right) .$

Since all the conditions of Theorem \ref{t0} are satisfied, we obtain that
the function $f(z,t)$ is a Loewner chain or has an analytic and univalent
extension to the whole unit disk $\mathbb{U},$ for all $t\in \left[ {%
0,\infty }\right) .$

From (\ref{3.2})-(\ref{3.4}) it follows in particular that%
\begin{equation*}
f(z,0)=\frac{v(\frac{1}{z})}{u(\frac{1}{z})}=\frac{1}{f(\frac{1}{z})}\in 
\mathcal{S}
\end{equation*}%
\noindent and for $\frac{1}{z}=\zeta \in \mathbb{U}^{\ast }$ we\ say that $%
f(\zeta )$ is univalent in $\mathbb{U}^{\ast }$. Thus the proof is completed.
\end{proof}

For $\alpha =0$ in Theorem \ref{t1} we obtain following new result:

\begin{corollary}
\label{c1} Let \textbf{\ }$f\in \sum $ be locally univalent function in $%
\mathbb{U}^{\ast }.$ If there exists an analytic function $h$ with $\func{Re}%
h(\zeta )\geq \frac{1}{2}$ in $\mathbb{U}^{\ast }$ and $h(\zeta )=1+\frac{%
h_{2}}{\zeta ^{2}}+...$ such that%
\begin{equation}
\left\vert {\frac{1-h(\zeta )}{h(\zeta )}\left\vert \zeta \right\vert
^{2}-(\left\vert \zeta \right\vert ^{2}-1)\left[ {\frac{\zeta {h}^{\prime
}(\zeta )}{h(\zeta )}+\frac{\zeta {f}^{\prime \prime }(\zeta )}{{f}^{\prime
}(\zeta )}}\right] }\right\vert \leqslant 1  \label{3.11}
\end{equation}%
\noindent for all$\ \zeta \in \mathbb{U}^{\ast },$ then $f$ is univalent in $%
\mathbb{U}^{\ast }.$
\end{corollary}

For $\alpha =\frac{1}{2}$ in Theorem \ref{t1} we obtain univalence criterion
given by Miazga and Wesolowski \cite{MiWe}.

\begin{corollary}
\label{c2}\textbf{\ }Let \textbf{\ }$f,g\in \sum $ be locally univalent
functions in $\mathbb{U}^{\ast }.$ If there exists an analytic function $h$
with $\func{Re}h(\zeta )\geq \frac{1}{2}$ in $\mathbb{U}^{\ast }$ and $%
h(\zeta )=1+\frac{h_{2}}{\zeta ^{2}}+...$ such that%
\begin{eqnarray}
&&\left\vert {\frac{1-h(\zeta )}{h(\zeta )}\left\vert \zeta \right\vert
^{2}-(\left\vert \zeta \right\vert ^{2}-1)\left[ {\frac{\zeta {h}^{\prime
}(\zeta )}{h(\zeta )}+\frac{\zeta {g}^{\prime \prime }(\zeta )}{{g}^{\prime
}(\zeta )}}\right] }\right.  \label{3.12} \\
&&\left. {\ +}\frac{{1}}{2}{(\left\vert \zeta \right\vert ^{2}-1)^{2}\frac{%
\zeta }{\overline{\zeta }}h(\zeta )\left[ {\left( {S_{f}(\zeta )-S_{g}(\zeta
)}\right) }\right] }\right\vert \leqslant 1  \notag
\end{eqnarray}%
\noindent for all$\ \zeta \in \mathbb{U}^{\ast },$ then $f$ is univalent in $%
\mathbb{U}^{\ast }.$
\end{corollary}

For $h(\zeta )=1$ and $\alpha =\frac{1}{2}$ in Theorem \ref{t1} we obtain
sufficient condition of Epstein type \cite{Eps} on the exterior of the unit
disk obtained earlier by Wesolowski\ \cite{We}.

\begin{corollary}
\label{c3} Let \textbf{\ }$f,g\in \sum $ be locally univalent functions in $%
\mathbb{U}^{\ast }.$ If the following inequality%
\begin{equation}
\left\vert \frac{{1}}{2}{(\left\vert \zeta \right\vert ^{2}-1)^{2}\frac{%
\zeta }{\overline{\zeta }}\left[ {\left( {S_{f}(\zeta )-S_{g}(\zeta )}%
\right) }\right] -(\left\vert \zeta \right\vert ^{2}-1){\frac{\zeta {g}%
^{\prime \prime }(\zeta )}{{g}^{\prime }(\zeta )}}}\right\vert \leqslant 1
\label{3.13}
\end{equation}%
\noindent is satisfied for all$\ \zeta \in \mathbb{U}^{\ast },$ then $f$ is
univalent in $\mathbb{U}^{\ast }.$
\end{corollary}

For $f(\zeta )=g(\zeta ),$ $h(\zeta )=1$ and $\alpha =\frac{1}{2}$ in
Theorem \ref{t1} we obtain well-known Becker's univalence criterion\ \cite%
{Be} in $\mathbb{U}^{\ast }.$

\begin{corollary}
\label{c4}\textbf{\ }Let \textbf{\ }$f\in \sum $ be locally univalent
function in $\mathbb{U}^{\ast }.$ If \ the following inequality%
\begin{equation}
{(\left\vert \zeta \right\vert ^{2}-1)}\left\vert {{\frac{\zeta {f}^{\prime
\prime }(\zeta )}{{f}^{\prime }(\zeta )}}}\right\vert \leqslant 1
\label{3.14}
\end{equation}%
\noindent is satisfied for all$\ \zeta \in \mathbb{U}^{\ast },$ then $f$ is
univalent in $\mathbb{U}^{\ast }.$
\end{corollary}

For $g(\zeta )=\zeta ,$ $h(\zeta )=1$ and $\alpha =\frac{1}{2}$ in Theorem %
\ref{t1} we obtain Nehari type univalence criterion \cite{Ne} in $\mathbb{U}%
^{\ast }.$

\begin{corollary}
\label{c5} Let \textbf{\ }$f\in \sum $ be locally univalent function in $%
\mathbb{U}^{\ast }.$ If the following inequality%
\begin{equation}
\left\vert {{{S_{f}(\zeta )}}}\right\vert \leqslant \frac{{2}}{{(\left\vert
\zeta \right\vert ^{2}-1)^{2}}}  \label{3.15}
\end{equation}%
\noindent is satisfied for all$\ \zeta \in \mathbb{U}^{\ast },$ then $f$ is
univalent in $\mathbb{U}^{\ast }.$
\end{corollary}

\textbf{Acknowledgement. }The present investigation was supported by Atat%
\"{u}rk University Rectorship under BAP Project (The Scientific and Research
Project of Atat\"{u}rk University) Project No: 2010/28.

Atat\"{u}rk University, Faculty of Science, Department of Mathematics,
25240, Erzurum-TURKEY.

edeniz@atauni.edu.tr; horhan@atauni.edu.tr


\begin{thebibliography}{99}
\bibitem{Be} J. Becker, Lownersche differentialgleichung und
schlichtheitskriterien, Math. Ann. 202 (1973) 321-335.

\bibitem{DeOr} E. Deniz, H. Orhan, Univalence criterion for analytic
functions, Gen. Math. 17 (4) (2009) 211--220.

\bibitem{Eps} C. L. Epstein, Univalence criteria and surfaces in hyperbolic
space, J. Reine. Angew. Math. 380 (1987) 196-214.

\bibitem{KaSr} S. Kanas, H. M. Srivastava, Some criteria for univalence
related to Ruscheweyh and S\u{a}l\u{a}gean derivatives, Complex Variables
Theory Appl. 38 (3) (1997) 263-275.

\bibitem{Le1} Z. Lewandowski, On a univalence criterion, Bull. Acad. Polon.
Sci. Math. 29 (1981) 123--126.

\bibitem{Le2} Z. Lewandowski, New remarks on a univalence criteria, Ann.
Univ. Mariae Curie-Sk\l odowska Sect. A Lublin-Polonia 41 (6) (1987) 41-49.

\bibitem{MiWe} J. Miazga, A. Wesolowski, A univalence criterion for
meromorphic functions, Ann. Polon. Math. 56 (1) (1991) 63-66.

\bibitem{Ne} Z. Nehari, The Schwarzian derivative and schlicht functions,
Bull. Amer. Math. Soc. 55 (1949) 545-551.

\bibitem{Po1} Ch. Pommerenke, Univalent functions\textit{,} Vandenhoech and
Ruprecht, G\"{o}ttingen (1975).

\bibitem{Po2} Ch. Pommerenke, \"{U}ber die Subordination analytischer
Funktionen, J. Reine~Angew. Math. 218 (1965) 159-173.

\bibitem{Ru} St. Ruscheweyh, An extension of Becker's univalence condition,
Math. Ann. 220 (1976) 285-290.

\bibitem{Sr} H. M. Srivastava, S. Owa (Editors), Current Topics in Analytic
Function Theory, World Scientifc Publishing Company, Singapore, New Jersey,
London and Hong Kong, 1992.

\bibitem{We} A. Wesolowski, On certain univalence criteria, Colloq. Math. 62
(1991) 39-42.
\end{thebibliography}
\end{document}